\documentclass[12pt]{amsart}

\usepackage{amssymb}
\newtheorem{lemma}{Lemma}
\newtheorem{example}{Example}
\newtheorem{remark}{Remark}
\newtheorem{definition}{Definition}
\newtheorem{theorem}{Theorem}
\newtheorem{corollary}{Corollary}
\newtheorem{proposition}{Proposition}

\begin{document}

{

\title[Varieties of chord diagrams]{On the cohomology of varieties of chord diagrams}
\author{V.A.~Vassiliev}
\address{Weizmann Institute of Science}
\email{vavassiliev@gmail.com}

\begin{abstract}
We study the space of codimension two subalgebras in $C^\infty(S^1, {\mathbb R})$ defined by pairs of conditions $f(\varphi)=f(\psi)$, $\varphi \neq \psi \in S^1$, or by their limits. We compute the mod 2 cohomology ring of this space, and also  the Stiefel--Whitney classes of the tautological vector bundle on it.
\end{abstract}

\thanks{This work was supported by the Absorption Center in Science of the Ministry of Immigration and Absorption of the State of Israel}
\maketitle
\unitlength 1mm

\section{Introduction}

\begin{definition} \rm
A {\em chord} is an unordered pair of points $\varphi \neq \psi \in S^1$.
A {\em chord diagram} is a finite collection of chords $\{\varphi_i, \psi_i\}$. The {\em rank} of a chord diagram is the codimension in $C^\infty(S^1, {\mathbb R})$ of the subalgebra consisting of the functions $f: S^1 \to {\mathbb R}$ which satisfy all conditions $f(\varphi_i) = f(\psi_i)$ over all chords of this diagram.  The space of all such subalgebras defined by chord diagrams of rank $n$ is denoted by $CD_n$. $\overline{CD_n}$ is the closure of $CD_n$ in the space of all subspaces of codimension $n$ in $C^\infty(S^1, {\mathbb R})$. The {\em canonical normal bundle} over $\overline{CD_n}$ is the $n$-dimensional vector bundle, whose fiber over a point is the quotient space of $C^\infty(S^1, {\mathbb R})$ by the corresponding subalgebra. $B(X, k)$ denotes the {\em configuration space} of all unordered collections of $k$ distinct points of the topological space $X$.
\end{definition}

Below we compute the cohomology ring of the space $\overline{CD}_2$ and the Stiefel--Whitney classes of the canonical normal bundle on it. For the motivation of this study, see \cite{pacific}, where in particular the Stiefel--Whitney classes of canonical bundles of spaces of chord diagrams in ${\mathbb R}^1$ were applied to problems in knot theory and interpolation theory. For some corollaries of our calculations, see Proposition \ref{corm} below.

\begin{example} \rm
\label{ex1}
The space $CD_1 \sim B(S^1, 2)$ is homeomorphic to the open Moebius band. 
Its closure $\overline{CD_1}$ is homeomorphic to the closed Moebius band and is obtained from $CD_1$ by adding all subalgebras parameterized by the points $\varphi \in S^1$ and consisting of functions satisfying the condition $f'(\varphi)=0$.
A cell decomposition of this space consists of four cells
\begin{picture}(8,6)
\put(4,2){\circle{6}}
\put(4,-0.9){\line(0,1){5.8}}
\put(7,2){\circle*{1}}
\end{picture}, 
\begin{picture}(8,6)
\put(4,2){\circle{6}}
\put(1,2){\line(1,0){6}}
\put(7,2){\circle*{1}}
\end{picture}, 
\begin{picture}(8,6)
\put(4,2){\circle{6}}
\put(0.3,0.8){\small$\ast$}
\put(7,2){\circle*{1}}
\end{picture}, 
\begin{picture}(8,6)
\put(4,2){\circle{6}}
\put(7,1){\small $\ast$}
\put(7,2){\circle*{1}}
\end{picture}
\ of dimensions 2, 1, 1, 0, defined respectively by chords not containing the distinguished point $\bullet \in S^1$, chords containing this point, the conditions $f'(\varphi)=0$ for any $\varphi \neq \bullet$, and the unique such condition for $\varphi=\bullet$. The mod 2 cellular complex with these generators has the boundary operators $\partial_2 \left( \begin{picture}(8,4)
\put(4,1){\circle{6}}
\put(4,-1.9){\line(0,1){5.8}}
\put(7,1){\circle*{1}}
\end{picture} \right) = \begin{picture}(8,4)
\put(4,1){\circle{6}}
\put(0.3,-0.2){\small$\ast$}
\put(7,1){\circle*{1}}
\end{picture}$, $\partial_1 \left( \begin{picture}(8,4)
\put(4,1){\circle{6}}
\put(1,1){\line(1,0){6}}
\put(7,1){\circle*{1}}
\end{picture} \right) = \partial_1 \left( \begin{picture}(8,4)
\put(4,1){\circle{6}}
\put(0.3,-0.2){\small$\ast$}
\put(7,1){\circle*{1}}
\end{picture}\right) = 0$, 
so its homology group $H_2$ is trivial, and $H_1 $ is generated by the class of the cell \begin{picture}(8,4)
\put(4,1.5){\circle{6}}
\put(1,1.5){\line(1,0){6}}
\put(7,1.5){\circle*{1}}
\end{picture} . Alternatively, the group $H_1(\overline{CD_1}, {\mathbb Z_2})$ is generated by the cycle 
\begin{picture}(8,4)
\put(4,1.5){\circle{6.3}}
\put(4,-1.5){\vector(0,1){6}}
\put(4,4.5){\vector(0,-1){6}}
\put(7,1.5){\circle*{1}}
\end{picture}
parametrized by the space $S^1/{\mathbb Z}_2 \equiv {\mathbb R}P^1$ of diameters of $S^1$: its elements are algebras of functions having equal values at the endpoints of the diameters.
\end{example}

\begin{example} \rm
The space $CD_2$ is naturally homeomorphic to the 2-configuration space $B(B(S^1,2), 2)$ of the 2-configuration space of $S^1$ (i.e. of the open Moebius band) factored by the following equivalence relation: for any three different points $\varphi, \psi, \chi \in S^1$, three pairs of points $((\varphi, \psi), (\varphi, \chi))$, $((\varphi, \psi), (\psi, \chi))$ and $((\varphi, \chi), (\psi, \chi))$ of $B(S^1,2)$ define the same point of $CD_2$.
\end{example}

The computations of \cite{pacific} imply that for any $n$ the Stiefel--Whitney classes $w_i$ of the canonical normal bundle on the space $CD_n$
are non-trivial for all $i \leq n - I_2(n)$, where $I_2(n)$ is the number of ones in the binary decomposition of $n$.

\begin{theorem}
\label{mthmadd}
The mod 2 cohomology groups of the space $\overline{CD_2}$
are isomorphic to ${\mathbb Z}_2$ in dimensions 0, 1, 2, and 3, and are trivial in all other dimensions. 
\end{theorem}

Let $W$ be the generator of the group $H^1(\overline{CD_2}, {\mathbb Z}_2)$.

\begin{theorem}
\label{mthmmult}
$W^2 \neq 0$ but $W^3 = 0$ in the ring $H^*(\overline{CD_2}, {\mathbb Z}_2)$.
\end{theorem}

\begin{corollary}
The group $H^i(\overline{CD_2}, {\mathbb Q})$ is isomorphic to ${\mathbb Q}$ for $i=0$ and $3$, and is trivial for all other $i$.
\end{corollary}

\noindent
{\it Proof.} The dimensions of rational cohomology groups are not greater than those with ${\mathbb Z}_2$ coefficients. By Theorem \ref{mthmmult}, the Bockstein operator of the generator of the group $H^1(\overline{CD_2}, {\mathbb Z}_2)$ is non-trivial, so it is a torsion group. The statement of the corollary then follows from the Euler characteristic  considerations. \hfill $\Box$

\begin{theorem}
\label{SWth}
1. The first Stiefel--Whitney class of the canonical normal bundle on $\overline{CD_1}$ is non-trivial.

2. The first Stiefel--Whitney class of the canonical normal bundle on $\overline{CD}_2$ is not trivial.

3. The second Stiefel--Whitney class of the canonical normal bundle on $\overline{CD}_2$ is non-trivial.
\end{theorem}

\begin{corollary}
\label{corSW}
The total Stiefel--Whitney class of the canonical normal bundle on $\overline{CD}_2$ is equal to $1+W+W^2$. \hfill $\Box$
\end{corollary}

Our calculations also imply the following Borsuk--Ulam-type statement.

\begin{proposition}
\label{BU}
For any pair of linearly independent smooth functions $f, g: S^1 \to {\mathbb R}$, there exists a non-trivial linear combination $\lambda f+ \mu g$ with real coefficients $\lambda, \mu$ whose derivative vanishes at a pair of opposite points of the circle. For a {\em generic} pair of functions $f$ and $g$, the number of such linear combinations $($considered up to multiplication by nonzero constants$)$ is odd.
\end{proposition}

For any pair of distinct chords, the space of all maps $S^1 \to {\mathbb R}^3$ taking equal values at the endpoints of each chord has codimension 6 in $C^\infty(S^1, {\mathbb R}^3)$; correspondingly, the intersection of this space with a generic 7-dimensional vector subspace $F^7 \subset C^\infty(S^1, {\mathbb R}^3)$ has dimension 1. 

\begin{proposition}
\label{corm}
For any 7-dimensional subspace $F^7 \subset C^\infty(S^1, {\mathbb R}^3)$, there exist pairs of distinct chords in $S^1$ such that the set of maps $f \in F^7$ that take equal values at the endpoints of each chord is at least two-dimensional. Moreover, the set of pairs of chords satisfying this condition is at least two-dimensional.
\end{proposition} 

\begin{remark} \rm The computations of \cite{pacific} imply only a weaker result: for any 6-dimensional subspace $F^6 \subset C^\infty(S^1, {\mathbb R}^3)$ there exists a pair of distinct chords in $S^1$ such that the set of maps $f \in F^6$ having equal values at the endpoints of each chord is at least one-dimensional.
\end{remark}

\section{Cell decomposition and homology group of $\overline{CD_2}$}
In the following pictures, each segment denotes the condition that the functions should take equal values at their endpoints; the asterisks denote the conditions of the vanishing derivative at the corresponding points. A double asterisk $\ast\ast$ at the point $\varphi$ denotes the condition $f'(\varphi)=f''(\varphi)=0$. In addition, the following subalgebras of codimension two appear in the variety $\overline{CD_2}$. 

For any ordered pair of points $\varphi \neq \psi \in S^1$ and a number $\alpha \in {\mathbb R}P^1$, the algebra $\gimel(\varphi, \psi; \alpha)$ consists of all functions $f$ such that $f(\varphi) = f(\psi)$ and $f'(\varphi) = \alpha f'(\psi)$. Obviously, $\gimel(\varphi, \psi; \alpha) \equiv \gimel(\psi, \varphi; \alpha^{-1})$. The three-dimensional cell $e^+$ (respectively, $e^-$) in $\overline{CD_2}$ consists of all such algebras with $\varphi \neq \bullet \neq \psi$ and $\alpha \in (0, +\infty)$ (respectively, $\alpha \in (-\infty, 0)$). The two-dimensional cells $e^+_\infty$ and $e^-_\infty$ are defined analogously, but with $\varphi = \bullet$. 

Also, for any point $\varphi \in S^1$ and number $\alpha \in {\mathbb R},$ the algebra $\circledast(\varphi; \alpha)$ consists of all functions such that $f'(\varphi)=0$ and $f'''(\varphi) = \alpha f''(\varphi)$. The two-dimensional cell $\Theta \subset \overline{CD}_2$ consists of all such algebras with $\varphi \neq \bullet$. The one-dimensional cell $\Theta_\infty$ consists of all such algebras with $\varphi=\bullet$.

\begin{proposition}
\label{mainprop}
The variety $\overline{CD_2}$ has the structure of a CW-complex with 
\begin{itemize}
\item three 4-dimensional cells $A, B, C$; 
\item nine 3-dimensional cells $a, b,$ $c,$ $d, e^+, e^-,$ $A_\infty, B_\infty, C_\infty;$
\item ten 2-dimensional cells $\Gamma$, $\Delta$, $\Xi$, $\Theta$, $a_\infty$, $b_\infty$, $c_\infty,$ $d_\infty$, $e^+_\infty$, $e^-_\infty$, 
\item five 1-dimensional cells $\aleph$, $\Gamma_\infty$, $\Delta_\infty$, $\Xi_\infty$, $\Theta_\infty$, 
\item and one 0-dimensional cell $\aleph_\infty$, 
\end{itemize}
which are 
described either in the following pictures or above in this section. $($The endpoints of the chords and the positions of the asterisks that do not coincide with the distinguished point $\bullet$ are the parameters of the corresponding cells$)$.

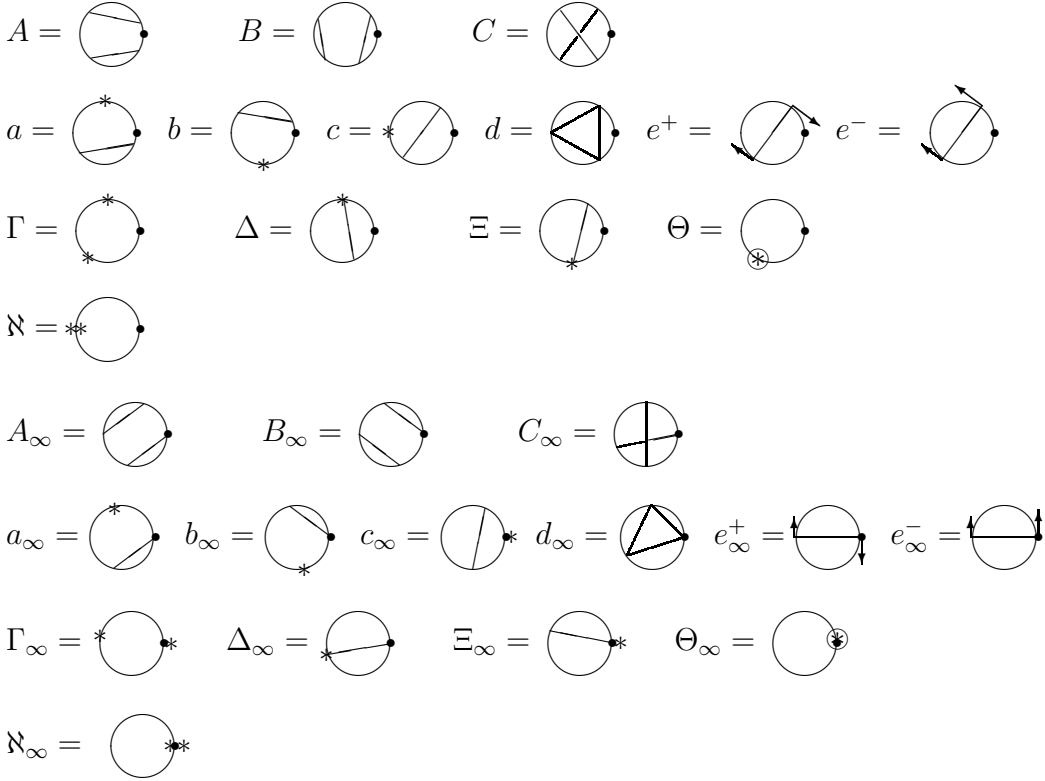
\unitlength 0.9mm
\noindent
$A= \mbox{\begin{picture}(12,10)
\put(6,1){\circle{10}} \put(10.8,1){\circle*{1}}
\put(3,-2.5){\line(6,1){7}}
\put(2.8,4.2){\line(5,-1){7.4}}
\end{picture}
}$ \qquad 
$B=\mbox{\begin{picture}(12,10)
\put(6,1){\circle{10}} \put(10.8,1){\circle*{1}}
\put(3,-2.7){\line(-1,6){1}}
\put(8,-3.2){\line(1,4){1.75}}
\end{picture}
}$ \qquad 
$C=\mbox{\begin{picture}(12,10)
\put(6,1){\circle{10}} \put(10.8,1){\circle*{1}}
\bezier{80}(3.3,-2.6)(4.5,-1)(5.7,0.6)
\bezier{80}(8.7,4.6)(7.5,3)(6.3,1.4)
\put(8.85,-2.75){\line(-3,4){5.5}}
\end{picture}
}$
\bigskip

\noindent
$a=\mbox{\begin{picture}(12,10)
\put(6,1){\circle{10}} 
\put(10.8,1){\circle*{1}}
\put(5,4.6){\small $\ast$}
\put(2.3,-2){\line(6,1){8}}
\end{picture}
}$ $b=\mbox{\begin{picture}(12,10)
\put(6,1){\circle{10}} \put(10.8,1){\circle*{1}}
\put(5,-4.8){\small $\ast$}
\put(2.3,4){\line(6,-1){8}}
\end{picture}
}$
$c=\mbox{\begin{picture}(12,10)
\put(6,1){\circle{10}} \put(10.8,1){\circle*{1}}
\put(0,0){\small $\ast$}
\put(3.3,-2.6){\line(3,4){5.5}}
\end{picture}
}$ 
$d=\mbox{\begin{picture}(12,10)
\put(6,1){\circle{10}} \put(10.8,1){\circle*{1}}
\bezier{120}(1.2,1)(4.8,3)(8.4,5)
\bezier{120}(1.2,1)(4.8,-1)(8.4,-3)
\bezier{120}(8.4,-3)(8.4,1)(8.4,5)
\end{picture}
}$ 
$e^+ =\mbox{
\begin{picture}(12,10)
\put(6,1){\circle{10}} \put(10.8,1){\circle*{1}}
\put(3,-3){\line(3,4){6}}
\put(3,-3){\vector(-4,3){3}}
\put(9,5){\vector(4,-3){4}}
\bezier{50}(3,-3)(1.5,-1.87)(0,-0.75)
\end{picture}
}$ $e^- =\mbox{
\begin{picture}(12,10)
\put(6,1){\circle{10}} \put(10.8,1){\circle*{1}}
\put(3,-3){\line(3,4){6}}
\put(3,-3){\vector(-4,3){3}}
\put(9,5){\vector(-4,3){4}}
\bezier{50}(3,-3)(1.5,-1.87)(0,-0.75)
\end{picture}
}
$
\bigskip

\noindent
$\Gamma= \mbox{\begin{picture}(12,10)
\put(6,1){\circle{10}} \put(10.8,1){\circle*{1}}
\put(2,-4){\small $\ast$}
\put(5,4.6){\small $\ast$}
\end{picture}
}$ \qquad 
$\Delta = \mbox{\begin{picture}(12,10)
\put(6,1){\circle{10}} \put(10.8,1){\circle*{1}}
\put(5,4.6){\small $\ast$}
\put(6.2,5.8){\line(1,-6){1.5}}
\end{picture}
}$ \qquad
$\Xi=\mbox{\begin{picture}(12,10)
\put(6,1){\circle{10}} \put(10.8,1){\circle*{1}}
\put(5,-5){\small $\ast$}
\put(6.2,-3.8){\line(1,4){2.2}}
\end{picture}
}$ \quad
$\Theta
=\mbox{\begin{picture}(12,10)
\put(6,1){\circle{10}} \put(10.8,1){\circle*{1}}
\put(2,-4.3){$\circledast$}
\end{picture}
}
$
\bigskip

\noindent
$\aleph = \mbox{\begin{picture}(12,10)
\put(6,1){\circle{10}} \put(10.8,1){\circle*{1}}
\put(-0.5,0){\small $\ast$}
\put(1,0){\small $\ast$}
\end{picture}
}$
\bigskip 

\smallskip

\noindent
$A_{\infty}=\mbox{\begin{picture}(12,10)
\put(6,1){\circle{10}} \put(10.8,1){\circle*{1}}
\put(11,1){\line(-4,-3){6.1}}
\put(1.2,1){\line(4,3){6}}
\end{picture}
}$ \qquad 
$B_{\infty}=\mbox{\begin{picture}(12,10)
\put(6,1){\circle{10}} \put(10.8,1){\circle*{1}}
\put(11,1){\line(-4,3){6.1}}
\put(1.2,1){\line(4,-3){6}}
\end{picture}
}$ \qquad $C_{\infty}= \mbox{\begin{picture}(12,10)
\put(6,1){\circle{10}} \put(10.8,1){\circle*{1}}
\put(11,1){\line(-5,-1){4.6}}
\bezier{80}(5.6,-0.1)(3.65,-0.5)(1.7,-0.9)
\put(6,-3.7){\line(0,1){9.5}}
\end{picture}
}$

\bigskip

\noindent
$a_{\infty}=\mbox{\begin{picture}(11,10)
\put(5,1){\circle{10}} \put(10,1){\circle*{1}}
\put(2.9,4.1){\small $\ast$}
\put(10,1){\line(-4,-3){6.1}}
\end{picture}
}$ $b_{\infty}=\mbox{\begin{picture}(11,10)
\put(5,1){\circle{10}} \put(10,1){\circle*{1}}
\put(5,-4.8){\small $\ast$}
\put(10,1){\line(-4,3){6.1}}
\end{picture}
}$ $c_{\infty}=\mbox{\begin{picture}(11,10)
\put(5,1){\circle{10}} \put(10,1){\circle*{1}}
\put(9.6,0){\small $\ast$}
\put(5,-3.7){\line(1,5){1.8}}
\end{picture}
}$ 
$d_{\infty}=\mbox{\begin{picture}(11,10)
\put(5.2,1){\circle{10}} \put(10,1){\circle*{1}}
\bezier{120}(9.8,1)(5.6,-0.35)(1.4,-1.7)
\bezier{80}(9.8,1)(7.4,3.4)(5,5.8)
\bezier{110}(5,5.8)(3.2,2.25)(1.4,-1.7)
\end{picture}
}$ $e^+_\infty
=\mbox{\begin{picture}(11,10)
\put(5,1){\circle{10}} \put(10,1){\circle*{1}}
\put(10,1){\line(-1,0){10}}
\put(10,1){\vector(0,-1){4}}
\put(0,1){\vector(0,1){3}}
\end{picture}
}
$ $e^-_{\infty}
=\mbox{\begin{picture}(10,10)
\put(5,1){\circle{10}} \put(10,1){\circle*{1}}
\put(10,1){\line(-1,0){10}}
\put(10,1){\vector(0,1){4}}
\put(0,1){\vector(0,1){3}}
\end{picture}
}
$
\bigskip

\noindent
$\Gamma_{\infty}= \mbox{\begin{picture}(12,10)
\put(6,1){\circle{10}} \put(10.8,1){\circle*{1}}
\put(10.8,-0.1){\small $\ast$}
\put(0.3,1){\small $\ast$}
\end{picture}
}$ \quad 
$\Delta_{\infty} = \mbox{\begin{picture}(12,10)
\put(6,1){\circle{10}} \put(10.8,1){\circle*{1}}
\put(0.2,-1.8){\small $\ast$}
\put(11,1){\line(-6,-1){10}}
\end{picture}
}$ \quad
$\Xi_{\infty}=\mbox{\begin{picture}(12,10)
\put(6,1){\circle{10}} \put(10.8,1){\circle*{1}}
\put(11,0){\small $\ast$}
\put(11,1){\line(-5,1){9.4}}
\end{picture}
}$ \quad
$\Theta_\infty =\mbox{\begin{picture}(12,10)
\put(6,1){\circle{10}} \put(10.8,1){\circle*{1}}
\put(9,0.4){$\circledast$}
\end{picture}
}$
\bigskip

\noindent
$\aleph_\infty = \mbox{
\begin{picture}(12,10)
\put(6,1){\circle{10}} \put(10.8,1){\circle*{1}}
\put(9,0){\small $\ast$}
\put(11,0){\small $\ast$}

\end{picture}
}$
\bigskip

The boundary operators of this cell complex mod 2 are as follows.

$\partial(A)=a + b + d + A_\infty + B_\infty,$ 

$\partial(B)=c + B_\infty + A_\infty + e^-, $ 

 $\partial(C)=d+e^+;$
\medskip

$\partial(a)=\Gamma+\Delta+a_\infty + c_\infty,$ 

$\partial(b)=c_\infty+\Xi + \Gamma + b_\infty,$

$\partial(c)=b_\infty+ \Xi + \Delta + a_\infty + \Theta,$ 

$\partial(d)=\Xi + \Delta,$ 

$\partial(e^+)=\Delta+\Xi,$ 

$\partial(e^-)=\Delta + \Xi + \Theta;$
\medskip

$\partial(\Gamma)=\aleph, $ 

$\partial(\Delta)=\Delta_\infty+\Xi_\infty+\aleph,$ 

$\partial(\Xi)=\Xi_\infty + \Delta_\infty + \aleph,$ 

$\partial(\Theta) = 0;$
\medskip 

$\partial (\aleph) = 0$
\bigskip

$\partial(A_\infty) = c_\infty + a_\infty + e^-_\infty, $ 

$\partial(B_\infty) = b_\infty + c_\infty+ e^-_\infty, $ 

$\partial(C_\infty) = 0;$
\medskip

$\partial(a_\infty) =\Gamma_\infty+ \Delta_\infty + \Xi_\infty + \Theta_\infty,$

$\partial(b_\infty) = \Xi_\infty+ \Delta_\infty + \Gamma_\infty+\Theta_\infty,$

$\partial(c_\infty) =\Gamma_\infty + \Theta_\infty,$

$\partial(d_\infty) =\Delta_\infty,$

$\partial(e^+_\infty) =\Delta_\infty + \Xi_\infty,$

$\partial(e^-_\infty) =\Delta_\infty + \Xi_\infty;$
\medskip

$\partial(\Gamma_\infty)= 0,$ 

$\partial(\Delta_\infty)=0,$

$\partial(\Xi_\infty)=0, $

$\partial(\Theta_\infty)=0.$
\end{proposition}

\noindent{\it Proof}: elementary calculations. \hfill $\Box$ \medskip

These formulas immediately imply the following detailing of Theorem \ref{mthmadd}.

\begin{corollary}
$H_4(\overline{CD_2}, {\mathbb Z}_2) \simeq 0$. The group $H_3(\overline{CD_2}, {\mathbb Z}_2)$ is isomorphic to $ {\mathbb Z}_2$ and is generated by the class of the cell $C_\infty$. The group $H_2(\overline{CD_2}, {\mathbb Z}_2)$ is isomorphic to ${\mathbb Z}_2$ and is generated by the chain $e^+_\infty + e^-_\infty$. The group $H_1(\overline{CD_2}, {\mathbb Z}_2)$ is isomorphic to ${\mathbb Z}_2$ and is generated by either of the cells $\Gamma_\infty$ or $\Theta_\infty$. $H_0(\overline{CD_2}, {\mathbb Z}_2) \simeq  {\mathbb Z}_2$. \hfill $\Box$ $\Box$
\end{corollary}

\section{Other realizations of homology groups}

Define the one-dimensional cycle $\tilde \Gamma_\infty \subset \overline{CD_2}$ parameterized by the projective line $S^1/{\mathbb Z_2}$ of pairs of opposite points of $S^1$: for each such pair we take the algebra consisting of functions whose derivative vanishes at these two points.

Define also the 2-dimensional cycle $\tilde e_\infty \subset \overline{CD_2}$ fibered over $S^1/{\mathbb Z}_2$, whose fiber over any pair of opposite points $(\varphi, \varphi+\pi) \subset S^1$ consists of all algebras $\gimel(\varphi, \varphi+\pi; \alpha)$, $\alpha \in {\mathbb R}P^1$. 
It is easy to see that this fiber bundle is non-orientable and thus is homeomorphic to the Klein bottle.

\begin{proposition}
\label{nz}
The first Stiefel--Whitney class of the canonical normal bundle on $\overline{CD_2}$ takes the non-zero value on the cycle $\tilde \Gamma_\infty$.
The second Stiefel--Whitney class of this bundle takes the non-zero value on the cycle $\tilde e_\infty \subset \overline{CD_2}$.
\end{proposition}

\begin{corollary}
\label{alter}
The cycle $\tilde \Gamma_\infty$ is homologous to the cycle $\Gamma_\infty$. The cycle $\tilde e_\infty$ is homologous to the cycle $e^+_\infty + e^-_\infty$. 
\end{corollary}

\noindent {\it Proof.} In both cases, the two classes being compared are non-trivial elements of a group isomorphic to ${\mathbb Z}_2$. \hfill $\Box$

\begin{proposition}
\label{regul}
 The non-trivial element of the group $H_2(\overline{CD_2}, {\mathbb Z}_2)$ can be represented by a two-dimensional cycle lying in $CD_2$. 
\end{proposition}

\noindent
{\it Proof.} Let $\varepsilon$ be a small positive number.
For any point $\gimel(\varphi, \varphi+\pi; \alpha) $ of the cycle $\tilde e_\infty$ consider the pair of chords $\left(\varphi +\varepsilon \frac{|\alpha|}{|\alpha|+1}, \varphi+\pi +\varepsilon \frac{1}{|\alpha|+1}\right)$ and $\left(\varphi -\varepsilon \frac{|\alpha|}{|\alpha|+1}, \varphi+\pi -\varepsilon \frac{1}{|\alpha|+1}\right)$ if $\alpha \in [0, +\infty]$, and the pair of chords $\left(\varphi +\varepsilon \frac{|\alpha|}{|\alpha|+1}, \varphi+\pi -\varepsilon \frac{1}{|\alpha|+1}\right)$ and $\left(\varphi -\varepsilon \frac{|\alpha|}{|\alpha|+1}, \varphi+\pi +\varepsilon \frac{1}{|\alpha|+1}\right)$ if $\alpha \in [-\infty, 0]$. 
These two chords never coincide (although they have a common endpoint if $\alpha=0$ or $\alpha=\infty$) and thus define a point of $CD_2$.
These formulas give the same result if we replace $\varphi$ by $\varphi+\pi$ and $\alpha$ by $\alpha^{-1}$, so they define a map $\tilde e_\infty \to CD_2$. This map is obviously homotopic to the identical embedding.  \hfill $\Box$

\section{Proof of Theorem \ref{SWth} and Propositions \ref{nz} and \ref{BU}.}

1. Ordering the endpoints $(\varphi, \varphi+\pi)$ of a diameter gives an orientation of the canonical normal bundle over the corresponding point of $CD_1$: the cosets of functions with $f(\varphi+\pi) > f(\varphi)$ belong to the positive part of it. Moving the point $\varphi$ continuously by the angle $\pi$ breaks this orientation.

2. Ordering the endpoints of such a diameter also specifies a canonical frame and thus an orientation of the canonical normal bundle over the corresponding point of the manifold $\tilde \Gamma_\infty \subset \overline{CD_2}$: its first (respectively, second) vector is the class of any function with $f'(\varphi+\pi)=1, f'(\varphi)=0$ (respectively,  $f'(\varphi+\pi)=0, f'(\varphi)=1$). Moving the point $\varphi$ by the angle $\pi$ permutes these two basic vectors and thus breaks the orientation.

3. Consider the section of the canonical normal bundle
over the manifold $\tilde e_\infty \subset \overline{CD_2}$ defined by the cosets of the function $\cos \varphi$. It has a single intersection point with the zero section of this bundle (i.e., a point $\gimel(\varphi, \varphi+\pi; \alpha) \in \tilde e_\infty$ containing this function). Indeed, the condition $\cos(\varphi+\pi)=\cos(\varphi)$ implies $\varphi = \pi/2 (\mbox{mod } \pi)$, and then necessarily $\alpha=-1$. This intersection is transversal, so the Euler characteristic of this bundle over $\tilde e_\infty$ is odd.  \hfill $\Box$
\medskip

\noindent 
{\bf Proof of Proposition \ref{BU}.} If a two-dimensional subspace of $C^\infty(S^1, {\mathbb R})$ does not contain non-trivial functions with derivative vanishing at two opposite points, then this subspace defines a trivialization of the restriction of the canonical normal bundle to the manifold 
$\tilde \Gamma_\infty \subset \overline{CD_2}$; this contradicts Proposition \ref{nz}.

\section{Proof of Theorem \ref{mthmmult}}

The surface $\tilde e_\infty \subset \overline{CD_2}$ is homeomorphic to the Klein bottle \unitlength 0.35 mm
\begin{picture}(30,30)
\put(0,0){\vector(1,0){30}}
\put(0,0){\vector(0,1){30}}
\put(0,30){\vector(1,0){30}}
\put(30,30){\vector(0,-1){30}}
\put(1,14){\footnotesize $u$}
\put(14,1){\footnotesize $v$}
\end{picture}
and is fibered over $S^1/{\mathbb Z}_2$. A generator $v$ of the group $H_1(\tilde e_\infty, {\mathbb Z}_2) \simeq {\mathbb Z}_2^2$ is given by a cross-section of this fiber bundle, namely it is swept out by the algebras $\gimel(\varphi, \varphi+\pi; 1)$, $\varphi \in [0,\pi)$. The space of functions 
\begin{equation} 
\label{fou1}
\lambda_1 \cos \varphi + \lambda_2 \sin \varphi \end{equation}
defines a trivialization of the canonical normal bundle over this generator: indeed, no non-trivial function of this form can satisfy both conditions $f(\varphi) = f(\varphi+\pi)$ and $f'(\varphi) = f'(\varphi+\pi)$ at some point $\varphi$. So the cohomology class $W$ takes zero value on this generator, which therefore represents a zero homology class of $\overline{CD_2}$.

Another generator $u$ of the group $H_1(\tilde e_\infty, {\mathbb Z}_2)$ 
consists of all algebras $\gimel(\pi, 0; \alpha)$, $\alpha \in {\mathbb R}P^1$. Consider the tautological morphism of the constant bundle consisting of the functions (\ref{fou1}) to the canonical normal bundle (taking each such function to its coset modulo the corresponding subalgebra). This morphism degenerates over a point $\{\alpha\}$ of our fiber if and only if a nonzero function of this type belongs to the subalgebra $\gimel(\pi, 0; \alpha)$, i.e. $\lambda_1 \cos \pi + \lambda_2 \sin \pi = \lambda_1 \cos 0 + \lambda_2 \sin 0$ and $(\lambda_1 \cos \pi + \lambda_2 \sin \pi)' = \alpha (\lambda_1 \cos 0 + \lambda_2 \sin 0)'$. These conditions imply $\lambda_1=0$ and $\alpha = \frac{\sin'(\pi)}{\sin'(0)} = -1$. This degeneration is of multiplicity 1 (i.e. the determinant of our morphism considered as a function of $\alpha$ has a simple root there), so the class $w_1(F^2)$ takes nonzero value on the cycle $u$.

So the ring homomorphism 
\begin{equation}
\label{restr}
H^*(\overline{CD_2}, {\mathbb Z}_2) \to H^*(\tilde e_\infty, {\mathbb Z}_2)\end{equation} defined by the inclusion $\tilde e_\infty \hookrightarrow \overline{CD_n}$
sends the class $W \in H^1(\overline{CD_2}, {\mathbb Z}_2)$ 
to the 1-cohomology class of $\tilde e_\infty$, which takes the value 1 on the generator $u$ and the value $0$ on the generator $v$. The square of this cohomology class is non-trivial in $H^2(\tilde e_\infty, {\mathbb Z}_2)$. By the functoriality of the cup product, it coincides with the image of the class $W^2$ under the map (\ref{restr}), hence $W^2 \neq 0$ in $H^2(\overline{CD_2}, {\mathbb Z}_2)$. 

This is sufficient to prove Corollary \ref{corSW}. Therefore, the total Stiefel--Whitney class of the third Cartesian power of the canonical normal bundle is equal to $(1+W+W^2)^3 \equiv 1+W+W^3$. To prove the second statement of Theorem \ref{mthmmult}, it remains to prove the following lemma.

\begin{lemma}
The third Stiefel--Whitney class of the third Cartesian power of the canonical normal bundle on 
$\overline{CD_2}$ is trivial.
\end{lemma}

\noindent
{\it Proof}. Consider the six-dimensional subspace in $(C^{\infty}(S^1, {\mathbb R}))^3 \equiv C^{\infty}(S^1, {\mathbb R}^3)$ consisting of 
maps $(f_1, f_2, f_3)$ whose three components $f_i$ are all homogeneous Fourier polynomials of degree 1, i.e. of the form (\ref{fou1}). This subspace defines a trivialization of the restriction of the third Cartesian power of the canonical normal bundle
to the closure of the cell  $C_\infty$, which generates the group $H_3(\overline{CD_2}, {\mathbb Z}_2)$: indeed, no non-trivial function of the form (\ref{fou1}) can belong to a subalgebra of the class $\overline{C_\infty}$.   \hfill $\Box \Box$

\section{Proof of Proposition \ref{corm}}

For any point of $\overline{CD_2}$, we have a natural homomorphism of the 
constant bundle with fiber $F^7$ to the third Cartesian product of the fiber of the canonical normal bundle at that point: the factorization modulo the space of all maps $f \equiv (f_1, f_2, f_3): S^1 \to {\mathbb R}^3$, where all three components $f_i$ belong to the corresponding subalgebra. The point of $CD_2$, i.e. a 2-chord diagram, satisfies the assertion of Proposition \ref{corm} if this map is not surjective or, equivalently, the conjugate map of dual vector spaces is not injective. If this never happens, then the third degree of the canonical normal bundle is isomorphic to a subbundle of the trivial 7-dimensional bundle, and has a one-dimensional orthogonal subbundle in it. The total Stiefel--Whitney class of this orthogonal subbundle is equal to the 
minus third degree of the total Stiefel-Whitney class of the canonical normal bundle over $CD_2$. By Theorems \ref{mthmmult} and \ref{SWth}, this is equal to the restriction of the class $(1 + W + W^2)^{-3} = (1+ W)^3 = 1 + W + W^2 \in H^*(\overline{CD_2}, {\mathbb Z}_2)$ to $CD_2$. By Proposition \ref{regul}, there exists a two-dimensional homology class of $CD_2$ such that the summand $W^2$ of this Stiefel--Whitney class takes non-zero value on any 2-cycle realizing this homology class. \hfill $\Box$

}

\end{document}